\documentclass[11pt]{amsart}
\usepackage{graphicx} % Required for inserting images
\usepackage[utf8]{inputenc}
\usepackage{amsmath}
\usepackage{xcolor}
\usepackage{amssymb}
\usepackage{amsmath,stackengine}
\usepackage{amsfonts}
\usepackage{amstext}
\usepackage{amsthm}
\usepackage[english]{babel}
\usepackage[autostyle]{csquotes}
\usepackage[super]{nth}
\usepackage{hyperref}
\usepackage[margin=1.37in]{geometry}
\usepackage{txfonts}

\newcommand{\W}{W_{0}^{1,p}}
\newcommand{\D}{\nabla}

\newcommand{\ul}{u_{\ell}}
\newcommand{\ui}{u_{\infty}}
\newcommand{\1}{\omega_{1}}
\newcommand{\2}{\omega_{2}}
\newcommand{\lamel}{\lambda_{\ell}^{1}}
\newcommand{\muinf}{\mu_{\infty}}
\newcommand{\norm}[1]{\left\lVert#1\right\rVert}

\newtheorem{Thm}{Theorem}[section]
\newtheorem{Lem}[Thm]{Lemma}

\newtheorem{example}[Thm]{Example}
\newtheorem{Def}[Thm]{Definition}

{%
\setcounter{enumi}{0}

\begin{enumerate}}%
{\end{enumerate} }  

\author[L. Esposito ]{Luca Esposito }
\address{\parbox{.8\linewidth}
{{\textbf{L. Esposito} (corresponding author)}\medskip \\
Dipartimento di Matematica, Università degli Studi di Salerno, Via Giovanni Paolo II 132, Fisciano 84084, Italy
 \medskip}}
\curraddr{}
\email{luesposi@unisa.it}

\author[L. Lamberti]{Lorenzo Lamberti}
\address{\parbox{.8\linewidth}
{{\textbf{L. Lamberti}}\medskip \\
Université de Lorraine, CNRS, IECL, F-54000 Nancy, France
 \medskip}}
\curraddr{}
\email{lorenzo.lamberti@univ-lorraine.fr}

\author[Dattatreya N. N]{Dattatreya N. N.}
\address{\parbox{.8\linewidth}
{{\textbf{Dattatreya N. N}}\medskip \\
Indian Institute of Technology - Kanpur, India \medskip}}
\curraddr{}
\email{dattatreya21@iitk.ac.in} 

\author[P. Roy]{Prosenjit Roy}
\address{\parbox{.8\linewidth}
{{\textbf{P. Roy}}\medskip \\
Indian Institute of Technology - Kanpur, India \medskip}}
\curraddr{}
\email{prosenjit@iitk.ac.in}

\makeatletter
\newcommand{\myitem}[1]{%
\medskip \item[#1]\protected@edef\@currentlabel{#1}%
} 

\hypersetup{
    colorlinks=true,
    linkcolor= blue,
    citecolor=blue,
    pdfpagemode=FullScreen,
    }

\def\om{\Omega}

\def\R{\mathbb R}

\def\XXint#1#2#3{{\setbox0=\hbox{$#1{#2#3}{\int}$}
     \vcenter{\hbox{$#2#3$}}\kern-.5\wd0}}

\title{Finsler $p$-Laplacian in domains becoming unbounded}

\begin{document}
\begin{abstract}
    We study the asymptotic behavior of sequences of solutions, energies functionals, and the first eigenvalues associated with the \textit{Finsler p-Laplace operator}, also known as the \textit{anisotropic p-Laplace operator} on a sequence of bounded cylinders whose length tends to infinity.
    We prove that the solutions on the bounded cylinders converge to the solution on the cross-section, with a polynomial rate of convergence in the general case and exponential convergence in some special cases. We show that energies on finite cylinders, with the multiplication of a scaling factor, converge to the energy on the cross-section. Finally, we investigate the convergence of the first eigenvalue and, for a specific subclass, we provide the optimal convergence rate.

\end{abstract}
\maketitle
\section{Introduction}
Over the past two decades, the behaviour of solutions defined on a sequence of bounded cylindrical domains that converge to an infinite cylinder has been extensively studied for a variety of equations, including elliptic, semi-linear elliptic, $p$-Laplace, pseudo-$p$-Laplace, anisotropic $p$-Laplace, and quasilinear equations, among others. 
Similar asymptotics analyses have also been carried out for the energy of the corresponding PDEs and the eigenvalue problems of various operators. This article extends these investigations to an anisotropic problem, the Finsler $p$-Laplacian arising by considering the Minkowski metric on $\R^{n}$, i. e.
\begin{equation*}
    \Delta_{H}^{p} u= -\mathrm{div}\left(H^{p-1}\left(\D u\right)\D H\left(\D u\right)\right).
\end{equation*}
Here `div' is the divergence operator and $H\colon\R^n\rightarrow[0,+\infty)$ is a function such that 

\begin{itemize}
\item $H(x)>0$ for $x\ne 0$.
     \item $H$ is convex.
    \item Homogeneity of degree-1: $H\left(tx\right)=|t|~H\left(x\right)$ for all $t\in \R$ and $x\in \R^{n}$.
    \item $H\in C^{1}\left(\R^n\setminus \{0\}\right)$.
\end{itemize}
\begin{example}($q-$ norms)
    The simplest examples are the $q$-norms on $\R^{n}$. For $1\leq q<+\infty$,
    \begin{equation*}
H=H_{q}\left(x\right):=\left(\displaystyle\sum_{i=1}^{n}|x_{i}|^{q}\right)^{\frac{1}{q}}.
    \end{equation*}
    From this, one can obtain 
    \begin{itemize}
      \item the Laplacian operator
       $$ -div\left(\D \cdot\right)=\Delta=\Delta_{H_{2}}^{2},$$
    \item the $p$-Laplacian operator
    $$-div\left(|\D\cdot|^{p-2}\D \cdot\right)=\Delta_{p}=\Delta_{H_{2}}^{p},$$
    \item the pseudo $p$-Laplacian operator $$\displaystyle\sum_{i=1}^{n}\frac{\partial~}{\partial x_{i}}\left(\left|\frac{\partial ~}{\partial x_{i}}\right|^{p-2}\frac{\partial ~}{\partial x_{i}}\right)=\Delta_{H_{p}}^{p}.$$
    \end{itemize}
    \end{example}
    
    \begin{example}
        Let $A$ be $n\times n$ matrix, then 
    \begin{equation*}
        H=H_{A,q}\left(x\right):=H_{q}\left(Ax\right),
    \end{equation*}
    where $H_{q}$ is defined in the previous example.
   One can get
   \begin{itemize}
       \item the constant coefficient elliptic operator
       $$- div\left(A\D u\right)=\Delta_{H_{A,2}}^{2}.$$
   \end{itemize}
    \end{example}
    \begin{example}
    Let $n=m_{1}+\cdots +m_{r}$ and $1\leq q,p_{i},\lambda_{i}<+\infty$ for $1\leq i\leq r$, denote $P=\left(p_{1},\cdots,p_{r}\right)$ and $\Lambda=\left(\lambda_{1},\cdots,\lambda_{r}\right)$. Consider 
    \begin{equation*}
        H=H_{q,P,\Lambda}:=\left(\displaystyle \sum_{i=1}^{r}\lambda_{i}\left(\sum_{j=s_{i-1}}^{m_{i}}|x_{j}|^{p_{i}}\right)^{\frac{q}{p_{i}}}\right)^{\frac{1}{q}},
    \end{equation*}
   % Then, one can get 
   % \begin{itemize}
       % \item Anisotropic p-Laplacian $$\displaystyle\sum_{i=1}^{r}\left(\lambda_{i}\sum_{j=1}^{m_{i}}\frac{\partial~}{\partial x_{j}}\left(\left|\frac{\partial ~}{\partial x_{j}}\right|^{p_{i}-2}\frac{\partial ~}{\partial x_{j}}\right)\right)=\Delta_{H_{q,P,\Lambda}}^{2}$$.
   % \end{itemize}
   where $s_{0}=0$ and $s_{k}=\displaystyle\sum_{i=1}^{k}m_{i}$.
\end{example}

However, the Finsler $p$-Laplacian does not generalize the $(p, q)$-Laplacian, which is the sum of the $p$-Laplacian and the $q$-Laplacian. In this case, one can consider $\Delta_{H_{2}}^{p}+\Delta_{H_{2}}^{q}$. Similarly, the Finsler $p$-Laplacian does not generalize the anisotropic $p$-Laplacian of the form
\begin{equation*}
    \displaystyle\sum _{i=1}^{n}\frac{\partial}{\partial x_{i}}\left(\left|\frac{\partial }{\partial x_{i}}\right|^{p_{i}-2}\frac{\partial}{\partial x_{i}}\right),
\end{equation*}
where $1\leq p_{i}<\infty$.\smallskip

Anisotropic problems emerge in a wide range of physical phenomenons, including minimal surface energy \cite{Taylor1978}, crystallography, crystal growth \cite{TaylorChanHandwerker1992}, image processing \cite{RudinOsherFatemi1992} and the Willmore functional \cite{Ulrich2004}. Additionally, nonlinear phenomena appear in fluid dynamics, such as porous media flow and flow of non-Newtonian fluid. To explore these problems, we consider the energy associated with the operator $\Delta_{H}^{p}$ on appropriate spaces, which is given by
 \begin{equation*}
     \int_{\om} H^{p}\left(\D u\left(X\right)\right) \ dX.
 \end{equation*}
Minimizing this energy when $H^{p}$ represents surface tension will produce the surface with the least energy \cite{Taylor1978}. A special case arises when surface tension is bounded below, known as \textit{Crystalline}, which is central to the study of the minimal surface for the surface energy $H$, also known as 
the \textit{Wulff Shape} \cite{DissertationXia}. Crystallines surface tension are those for which the Wulff Shape is polyhedral; in this direction, we refer to \cite{Wulff1901, ValdinociCozziAlberto2014, Giga2002} and the references therein. When the energy is constant in all directions, the problem reduces to the classical isoperimetric problem, as exemplified by soap bubbles. Further insight on the isoperimetric problem and Wulff shape can be found in \cite{BelloniFeronKawohl2003, Esposito2005} and references therein, and an overview on the genesis of anisotropic surface energy is available in \cite{EliJhon1977}. The same energy functional also generates the gradient flow of the crystalline energy. For a detailed explanation of this and the role of Wulff shape in crystal growth, we refer to \cite{NovagaPaolini1999} and the reference therein.\smallskip

In the mathematical framework, the Finsler Laplacian arises naturally in Finsler geometry. Consider a smooth manifold $M$ and its tangent bundle $TM$. A Finsler metric is given by a function $F(x,y)$ defined on $TM$ satisfying the properties mentioned previously with respect to $y$, with convexity replaced by \textit{strong convexity} of $F^{2}(x,\cdot)$, see \cite[Definition 2.2]{MezeiVas2019}. A manifold equipped with the Finsler metric is called a Finsler manifold \cite[Section 1.1]{BaoChernShen2000}. For a fixed $x\in M$, $F(x,\cdot)=H(\cdot)$ defines a \textit{Minkowski Norm} \cite[Section 1.2]{Shen2001} on the tangent space $T_{x}M$. However, $H$ need not be a norm (for an example of $H$, see \cite[Appendix A]{ValdinociCozziAlberto2014}). When $M=\Omega$, an open subset of $\R^{n}$, $T_{x}\Omega$ is identified with $\R^{n}$, thus the function $H$ in question is a Minkowski norm on $\R^{n}$, the Laplacian operator in this setting is $\Delta_{H_{2}}^{2}$.\smallskip

For $f\in L^{p/p-1}(\omega_2)$,  consider the Dirichlet problem 
\begin{equation}\label{eqn dirichlet problem on omel}
\begin{cases}
-\mathrm{div}\left(H^{p-1}\left(\D u\right)\D H\left(\D u\right)\right)=f(X_{2})\hspace{3mm} \text{in} \ \Omega_{\ell},\\
u\in W^{1,p}_0\left(\Omega_{\ell}\right).
\end{cases}
\end{equation}
We investigate the asymptotics behavior, as $\ell\to+\infty$ of the sequence of solutions $\{\ul\}_{\ell}$ and energies $J_{\ell}$ associated with \eqref{eqn dirichlet problem on omel}, posed on finite cylindrical domains $\om_{\ell}=\ell \1\times\2$. Here, each point is written as $x=(X_{1},X_{2})\in\om_{\ell}$, with $X_{1}\in \ell\1$ and $X_{2}\in \2$. We refer to Subsection \ref{notations} for more details about notation. 
Next, we turn to the variational eigenvalue problem 
\begin{equation}\label{equation eigenvalue problem}
    \begin{cases}
-\mathrm{div}\left(H^{p-1}\left(\D v\right)\D H\left(\D v\right)\right)=\lambda |v|^{p-2} v\hspace{3mm} \text{in}\ \Omega_{\ell},\\
v\in W^{1,p}_{0}\left(\Omega_{\ell}\right),
\end{cases}
\end{equation}
  and analyze the asymptotics behavior of the sequence of first eigenvalues $\left\{\lambda_{\ell}\right\}_{\ell}$ as $\ell\to+\infty$.\smallskip

 The study of equations on cylinders arises in various physical contexts, such as cylindrical gravitational waves generated by cylindrical sources propagating along the axial direction \cite{BuzzMatKir2023}, porous media flow \cite{BenHua2015}, particle accelerators and more. The behaviour of solutions defined on varying domains with a fixed cross-section has been investigated for different operators in various settings. For example, in the case of the elliptic operators including $\Delta^{2}_{H_{A,2}}$ with Dirichlet boundary conditions, the convergence rate of solutions on finite cylinders to the solution on the cross-section is shown to be polynomial in \cite[Theorem 1.1]{ChipotRougirel2002}. In addition, in \cite[Theorem 1.1]{chipotYeressian2008}, for the Laplacian $\Delta_{H_{2}}^{2}$, convergence is shown to be exponential. A semilinear elliptic equation involving $\Delta_{H_{2}}^{2}$ is addressed in \cite{ChipoJunaPino2017}, the eigenvalue problem for an elliptic equation including $\Delta^{2}_{H_{2}}$ is examined in \cite{ChipotRougirel2008, ChioptRoyShafi2013} and $\Delta^{p}_{H_{2}}$ in \cite{LucaRoyFiroj2021}. Studies related to $p$-Laplacian $\Delta^{p}_{H_{2}}$, anisotropic, and pseudo-$p$-Laplacian $\Delta_{H_{p}}^{p}$ for Dirichlet problems and eigenvalue problems can be found in \cite{ChipoXie2004, PJana2024}. In \cite[Theorem 1.1]{PJana2024} the convergence of the solution for pseudo-$p$-Laplace operator $\Delta_{H_{p}}^{p}$ is shown to be exponential, while the energy and the first eigenvalue convergence is polynomial. For problems involving infinite boundary data with $\Delta_{H_{p}}^{p}$, we refer to \cite{IndroDatta2024} and references therein. Additionally, we refer \cite{chipotYeressian2013, Mugnai2013, ChipotRoyMojsic2016} for energy-related convergence and variational inequality. In \cite{ChipotRoyMojsic2016}, the authors have considered the energy $\int_{\Omega_{\ell}}F(\nabla \ul)-f\ul\ dx$ where $F$ is a convex function of $q$-type with out smoothness assumption.\smallskip

 For the purposes of this article, we impose the following assumptions on $H$:
\begin{enumerate}
    \myitem{$\textbf{A1}$}\label{a1} \textit{For $p\geq 2$, 
\begin{equation*}
\label{monotonicity}
\left<H^{p-1}\left(z_1\right)\D H\left(z_1\right)-H^{p-1}\left(z_2\right)\D H\left(z_2\right)\cdot z_1-z_2\right >\geq c_1|z_1-z_2|^p.
\end{equation*}
For $1<p<2$,
\begin{align*}
\label{eq7}
&\left<H^{p-1}\left(z_1\right)\D H\left(z_1\right)-H^{p-1}\left(z_2\right)\D H\left(z_2\right) \cdot z_1-z_2\right>\geq c_2|z_1-z_2|^2\left(|z_1|+|z_2|\right)^{p-2},
%&\textcolor{blue}{\left<H^{p-1}\left(z_1\right)\D H\left(z_1\right)-H^{p-1}\left(z_2\right)\D H\left(z_2\right)\cdot z_1-z_2\right> \leq c_3|z_1-z_2|^2\left(|z_1|+|z_2|\right)^{p-2}}
\end{align*}
for any $z_1,z_2\in\R^n$ and for some positive constants $c_{1},c_2$.}
\myitem{$\textbf{A2}$}\label{a2} For $1< p<2$, there exists a positive constant $c_{3}$ such that for any any $z_{1}$ and $z_{2}$ in $\R^{n}$,
\begin{equation*}
    \left|H^{p-1}\left(z_1\right)\D H\left(z_1\right)-H^{p-1}\left(z_2\right)\D H\left(z_2\right)\right|\leq c_{3}|z_1-z_2|\left(|z_1|+|z_2|\right)^{p-2}.
\end{equation*}
%\myitem{$\textbf{A3}$}\label{a3} $$\displaystyle a:=\min_{|x|=1}H\left(x\right)>0$$
\end{enumerate}
Under these assumptions, we establish that the rate of convergence of $u_\ell$ is generally polynomial.
\begin{Thm}\label{thm poly convergence}
Let $\ui$ be the solution of 
\begin{equation}\label{eqn dirchlet prob cross-sec}
\begin{cases}
-\mathrm{div}\left(H^{p-1}\left(\D_{2} u_\infty\right)\D_{Z_{2}} H\left(\D_{2} u_\infty\right)\right)=f\left(X_{2}\right)\hspace{3mm} \text{in}\ \omega_2,\\
u_\infty\in W^{1,p}_0\left(\omega_2\right).
\end{cases}
\end{equation} For $p\geq 2$, it holds that
\begin{equation}
\norm{\D\left(u_\ell-u_\infty\right)}_{L^p\left(\Omega_{\ell/2}\right)}\leq C_{1}\ell^{m-\frac{p}{p-1}}.
\end{equation}
And For $1<p<2$, it holds that
\begin{equation}\norm{\D\left(u_\ell-u_\infty\right)}_{L^p\left(\Omega_{\ell/2}\right)}\leq C_{2}\ell^{m-\frac{p^2}{2-p}},
\end{equation}
where, $C_{i}$, for $~i=1,2$, is a positive constant and $u_{\infty}$, originally defined on $\omega_2$, is to be understood as extended constantly with respect to the variables $X_1$.
\end{Thm}

As noted earlier, the literature establishes that for $H=H_{p}$ with $p \geq 2$, the rate of convergence is exponential. In this regard, in Theorem \ref{thm expo convergence} (see Section \ref{notation and Preliminaries}), we provide sufficient conditions on 
$H$ that guarantee the exponential rate of convergence. \smallskip

Next, we study the asymptotic behaviour as $\ell$ goes to $\infty$ of the energy-functional 
\begin{equation*}
    J_{\ell}\left(v\right)=\int_{\om_{\ell}} \left\{\frac{H^{p}\left(\D v\right)}{p}-f\left(X_{2}\right)v \right\}dx,  \quad \text{for all } v\in W^{1,p}_{0}\left(\om_{\ell}\right),
\end{equation*}
of the problem \eqref{eqn dirichlet problem on omel}. In this respect, let 
\begin{equation*}
    J_{\infty}\left(w\right)=\int_{\2}\left\{\frac{H^{p}\left(\D_{X_2} w\right)}{p}-f\left(X_{2}\right)w\right\} \ dX_{2}, \quad \text{for all } w\in W^{1,p}_{0}\left(\2\right).
\end{equation*}
The main result in this direction is as follows.
\begin{Thm}\label{thrm energy}
    For any $1\leq p<\infty$, we have 
    \begin{equation*}
        \displaystyle\lim_{\ell\to\infty}\frac{J_{\ell}\left(u_{\ell}\right)}{\left|\ell \1\right|}=J_{\infty}\left(u_{\infty}\right).
    \end{equation*}
    Moreover, 
    \begin{equation*}
        J_{\infty}\left(u_{\infty}\right)\leq \frac{J_{\ell}\left(u_{\ell}\right)}{\left|\ell\1\right|}\leq J_{\infty}\left(u_{\infty}\right)+\frac{C}{\ell},
    \end{equation*}
    for some positive constant $C$ independent of $\ell$.
\end{Thm}

We conclude by studying the equation \eqref{equation eigenvalue problem}. 
Let $u_{\ell}$ and $u_{\infty}$ denote the positive minimizers of the following Rayleigh quotients:
\begin{equation}\label{raliegh quotient for lamda}
    \lambda_{\ell}^{1}=\min_{\substack{v \in W_0^{1,p}(\Omega_{\ell}), \ v \neq 0}}\frac{\int_{\om_{\ell}}H^{p}\left(\D v\right) \ dx}{\int_{\om_{\ell}}v^{p}\ dx}
\end{equation}

and 
\begin{equation}\label{raleigh quotient for mu}
    \mu_{\infty}= \min_{\substack{w \in W_0^{1,p}(\omega_2), \ w \neq 0}}\frac{\int_{\omega_{2}}H^{p}\left(\D_{2}w\right)\ dX_{2}}{\int_{\omega_{2}}w^{p}\ \ dX_{2}},
\end{equation}
respectively. Then, we derive the following result.
\begin{Thm}\label{eigenvalue from above}
    There exists a constant $C=C\left(p,\mu_{\infty}\right)$ independent of $\ell$ such that we have 
    \begin{equation*}
        \mu_{\infty}\leq\lambda_{\ell}^{1}\leq  \mu_{\infty}+\frac{C}{\ell}, 
    \end{equation*}
    for $\ell$ sufficiently large.
\end{Thm}
Further, for a subclass $\tilde{H}_{q}$ which includes the Laplacian and the pseudo-$p$-Laplacian, we show that the rate of convergence of $\lambda_{\ell}^{1}$ to $\mu_{\infty}$ is $\frac{1}{\ell^{p}}$, which is done in  Theorem \ref{eigenvalue lower bound withspeed}. To aid in the last result, we state and prove a version of Picone's identity  in Lemma \ref{picon}.\smallskip

The novelty of this article lies in proving new results for the Finsler Laplacian, which, as a consequence, allow for a unified treatment of several known results established over the past two decades for other operators such as \textit{Laplacian, constant coefficient elliptic operators, $p$-Laplacian and pseudo-$p$-Laplacian} that has been accomplished over the last two decades.\smallskip

The paper is structured as follows: we give notations, definitions and preliminary results in Section \ref{notation and Preliminaries}. In Section \ref{Asymptotics of solutions}, we prove Theorem \ref{thm poly convergence} and a sufficient condition for exponential convergence of the solutions $\ul$. Section \ref{Energy asymptotics} presents the proof of Theorem \ref{thrm energy} and Section \ref{An eigenvalue problem} contains the proof of Theorem \ref{eigenvalue from above} and the case of the exact rate of convergence of $\lambda_{\ell}^{1}$ to $\mu_{\infty}$. In section \ref{Appendix}, we briefly give the direct method of calculus of variation for the existence of a solution to \eqref{eqn dirichlet problem on omel} followed by the uniqueness.

\section{Notation and Preliminaries}\label{notation and Preliminaries}
We fix some notation for the rest of the article, define the notion of a solution and explain some basic results that will be useful.
\subsection{Notation}\label{notations}
\begin{itemize}
    \item We denote $\om_{\ell}=\ell \1\times\2$ for $\ell>0$, where, $\omega_1\subset \R^m$, $1\leq m<n$, is open, bounded and convex, with $0\in\omega_1$ and $\omega_2\subset \R^{n-m}$ is bounded and open.
    \item For any $x\in\R^n$ we write $x=\left(X_{1},X_{2}\right)$, where $X_1=\left(x_1,\cdots,x_m\right)$ and $X_2=\left(x_{m+1},\cdots,x_n\right)$.
    \item By $\D_{X_{i}}$, $i=1,2$ we mean the gradient with respect to the variables of $X_{i}$. And we define $\D_{1}:=\left(\partial_{x_{1}},\cdots,\partial_{x_{m}},0,\cdots,0\right)\in \R^{n}$ and $\D_{2}:=\left(0,\cdots,0,\partial_{x_{m+1}},\cdots,\partial_{x_{n}}\right)\in \R^{n}$.
    \item For any set $\Omega\in \R^{r}$, $\left|\om\right|$ denote the Lebesgue measure of $\om$ in the corresponding dimension $r$.
    \item For $p>1$ we denote by $p'$ the conjugate of $p$, that is $\frac{1}{p}+\frac{1}{p'}=1$.
    \item  We denote 
    \begin{equation}\label{b}
        \theta_{1}=\min_{|x|=1}H(x), \quad\quad\text{and} \quad\quad \theta_{2}=\displaystyle\max_{|x|=1}H\left(x\right).
    \end{equation}
    \item We denote by $C$ a generic positive constant independent of $\ell$ which may vary from line to line.
\end{itemize}

   We will use the following properties of $H$ in various stages of this article.
   
\begin{Lem}\label{lemma_properties of H}
    Let $H:\R^{n}\to [0,+\infty)$ be a convex, positively 1-homogenous and $C^{1}(\R^{n}\setminus\{0\})$ function. Then the following holds:
    \begin{enumerate}
        \item For $\theta_{1}$ and $\theta_{2}$ given in \eqref{b}, we have $\theta_{1}|z|\leq H\left(z\right)\leq \theta_{2}|z|$, for any $z\in\R^n$, and $0<\theta_{1}\leq \theta_{2}$.
        \item $H\left(x+y\right)\leq H\left(x\right)+H\left(y\right)$, for any $x,y\in\R^n$.
\item There exists a constant $c>0$ such that $|\D H|\leq c$.
\item For any $x\in\R^n$, $\left<\D H(x)\cdot x\right>=H\left(x\right)$.
\end{enumerate}
\end{Lem}
\begin{proof}
    From the homogeneity,  $H\left(z\right)=|z|H\left(\frac{z}{|z|}\right)$, the assertion (1) follows by comparing it with the maximum and minimum of $H$ on the unit circle. Strict inequality follows the continuity and strict positivity and homogeneity of $H$. When $\theta_{1}=\theta_{2}$, $H=\theta_{1}H_{2}$.
\smallskip

    The assertion (2) holds by taking $t=\frac{1}{2}$ in the definition of convexity and applying the homogeneity of $H$.\\
    \indent Now we prove the assertion (3). By $(2)$ and the homogeneity we have that
    \begin{equation*}
        \frac{H\left(x+he_{i}\right)-H\left(x\right)}{h}\leq H\left(e_{i}\right),\quad\forall h\neq 0. 
    \end{equation*}
    Thus the assertion (3) holds with $c=\left(\displaystyle\sum_{i=1}^{n}H^{2}(e_{i})\right)^{\frac{1}{2}}$.\\
    \indent The assertion (4) can be obtained by differentiating $H\left(tx\right)=t H\left(x\right), t>0$ with respect to $t$ and substituting $t=1$.
\end{proof}
The dual $H_{0}$ of $H$ is defined as
\begin{equation*}
    H_{0}\left(\xi\right):=\displaystyle\sup_{x\neq 0}\frac{\left<\xi\cdot x\right>}{H\left(x\right)} \quad \text{for all } \xi\in \R^{n}.
\end{equation*}
\begin{Lem}\label{H0} $H$ and $H_{0}$ are related through the following properties:
    \begin{enumerate}
     \item $H$-H\"older: $\left<\xi\cdot x\right>\leq H_{0}\left(\xi\right)H\left(x\right)$, for any $x,\xi\in\R^n$.
     \item $H_{0}$ is convex and positively 1-homogeneous.
        \item The dual of $H_{0}$ is $H$.
        \item $H_{0}\left(\D H\right)=1$ and $H\left(\D H_{0}\right)=1$.
    \end{enumerate}
\end{Lem}
\begin{proof}
    If $x=0$, then $H$-H\"older holds trivially, otherwise it holds by 
    \begin{equation*}
        \left<\xi \cdot x\right>=\frac{\left<\xi \cdot x\right>H\left(x\right)}{H\left(x\right)},
    \end{equation*}
    and the definition of $H_{0}.$
    \smallskip
    
    The second can be proved by using the definition.
    \smallskip
    
    For the proof of others, we refer to \cite[Proposition 1.3]{DissertationXia}, \cite[Section 2.1]{BePo1996} and the reference therein.
\end{proof}
\begin{Def}
    Let $f\in L^{p'}\left(\om\right)$, a function $u\in \W\left(\om\right)$ is a weak solution to 
    \begin{equation*}
\begin{cases}
-\mathrm{div}\left(H^{p-1}\left(\D u\right)\D H\left(\D u\right)\right)=f\left(X_{2}\right)\hspace{2mm} \text{in} \ \Omega,\\
u\in W^{1,p}_0\left(\Omega\right),
\end{cases}
\end{equation*}
    if 
    \begin{equation*}
        \int_{\om}H^{p-1}\left(\D u\left(x\right)\right)\left<\D H\left(\D u\left(x\right)\right)\cdot\D \phi\left(x\right)\right>\ dx =\int_{\om} f\left(X_{2}\right) \phi\left(x\right) \ dx, \quad \forall\phi\in \W\left(\om\right).
    \end{equation*}
\end{Def}

The existence of the solution follows from variational methods. For completeness, we give the idea of the proof of existence and uniqueness in the appendix (Section    \ref{Appendix}).

\begin{Def}
    A function $v\in W_{0}^{1,p}\left(\om\right)$ is a weak solution to 
    \begin{equation*}
    \begin{cases}
-\mathrm{div}\left(H^{p-1}\left(\D v\right)\D H\left(\D v\right)\right)=\lambda |v|^{p-2} v\hspace{0.1mm} &\text{in}\ \Omega_{\ell},\\
v\in W^{1,p}_{0}\left(\Omega_{\ell}\right),
\end{cases}
\end{equation*}if
    \begin{equation*}
        \int_{\om}H^{p-1}\left(\D v\right)\left<\D H\left(\D v\right)\cdot \D\phi\right>\ dx=\lambda \int_{\om}\left|v\right|^{p-2}v\phi \ dx, \quad \text{for all } \phi\in W_{0}^{1,p}\left(\om\right).
    \end{equation*}
\end{Def}
Then the positive minimizers $\ul$ and $\ui$ of \eqref{raliegh quotient for lamda}
and \eqref{raleigh quotient for mu}
satisfies the above weak formulation on the respective domains in the respective dimensions (see \cite[Theorem 3.1]{BelloniFeronKawohl2003}). \smallskip

Next, we state an auxiliary lemma.
\begin{Lem}\label{inequlity lemma}
    For $a,b\geq 0$ 
  \begin{enumerate}
  \item For $\alpha\geq 1$, we have 
    \begin{equation*}
        (a+b)^{\alpha}\leq 2^{\alpha-1}(a^{\alpha}+b^{\alpha}).
    \end{equation*}
      \item{(\cite[Lemma 2.5]{AdiRoyVikek})} Let $\alpha>1$. $a, b\in \R$ be non-negative. Then for every $c>1$ it holds that
      \begin{equation*}
        (a+b)^{\alpha}\leq ca^{\alpha}+ (1-c^{\frac{-1}{\alpha-1}})^{1-\alpha}b^{\alpha}.
    \end{equation*}
    \item For $\alpha\geq 1$, we have
    \begin{equation*}
        (a+b)^{\alpha}\geq a^{\alpha}+b^{\alpha}.
    \end{equation*}
    \item For $0\leq \alpha \leq 1$, we have 
    \begin{equation*}
        (a+b)^{\alpha}\leq a^{\alpha}+b^{\alpha}.
    \end{equation*}
  \end{enumerate}
    
\end{Lem} 
We provide a brief proof of \textit{Picone's identity}. The proof is straightforward and relies on Lemma \ref{H0}.\\
\begin{Lem}[{\cite[Theorem 3.1]{Jarovs2014}} Picone's identity]\label{picon}
    For any differentiable functions $u\geq 0,v>0$, we set
    \begin{equation*}
        R\left(u,v\right)= H^{p}\left(\D u\right)- \left<{H^{p-1}\left(\D v\right)\D H\left(\D v\right)}\cdot \D \left(\frac{u^{p}}{v^{p-1}}\right)\right>,
    \end{equation*}
    and 
    \begin{equation*}
        L\left(u,v\right)=H^{p}\left(\D u\right)+\left(p-1\right)H^{p}\left(\frac{u}{v}\D v\right)-p \left<H^{p-1}\left(\frac{u}{v}\D v\right)\D H\left(\D v\right) \cdot \D u\right>.
    \end{equation*}
    Then we have $R(u,v)=L(u,v)\geq 0$.
\end{Lem}
\begin{proof}
    Using (4) from Lemma \ref{lemma_properties of H}, by a simple expansion we get that $R(u,v)=L(u,v)$. Then by $H$-H\"older inequality, the homogeneity of $H_{0}$ and the action of $H_{0}$ on $\D H$, one has 
    \begin{equation*}
        \left< H^{p-1}\left(\frac{u}{v}\D v\right)\D H\left(\D v\right) \cdot \D u\right> \leq H_{0}\left(H^{p-1}\left(\frac{u}{v}\D v\right)\D H\left(\D v\right)\right)H\left(\D u\right)=H^{p-1}\left(\frac{u}{v}\D v\right) H\left(\D u\right).
    \end{equation*}
    The sign of $L(u,v)$ follows from Young's inequality.
\end{proof}
In \cite[Lemma 1.1]{Jarovs2014}, the author has proved a stronger version of the previous lemma in which $u,v$ are differentiable such that $v\neq 0$, in which case the corresponding $R$ need not coincide with $L$. We proved a weaker version that is sufficient to establish the lower bound in the Theorem \ref{eigenvalue from above}.\smallskip

Finally, we state the \textit{Poincar\'e inequality on strips}. We only outline the proof, which is straightforward.
\begin{Lem}[{\cite[Lemma 2.2]{ChipoXie2004}} Poincar\'e inequality on strips]\label{poincare on strip}
    There exists a constant $C>0$, not depending on $\ell$ such that 
    \begin{equation*}
        \norm{u}_{L^{p}(\om_{\ell})}\leq C \norm{\D_{X_{2}}u}_{L^{p}(\om_{\ell})},
    \end{equation*}
    for all $u\in W^{1,p}(\om_{\ell})$ with $u=0$ on $\ell\1\times\partial \2 $. 
\end{Lem}
\begin{proof}
By the Poincar\'e inequality on $\2$ we have
\begin{equation*}
    \int_{\2}|u\left(X_{1},\cdot\right)|^{p}\ dX_{2}\leq C(p,\2)\int_{\2}|\D_{X_{2}}u\left(X_{1},\cdot\right)|^{p}\ dX_{2}.
\end{equation*}
    The result follows by integrating over $\ell\1$.
\end{proof} 

\section{Asymptotics of solutions}\label{Asymptotics of solutions}
In this section, we give the proof of Theorem \ref{thm poly convergence} and establish the exponential rate of convergence for a subclass of $H$. We break down parts of the proof into the following lemmas where we provide some useful estimates. In the following lemmata $u_l$ denote a solution of the Dirichlet problem \eqref{eqn dirichlet problem on omel} while $u_{\infty}$ denote a solution of the Dirichlet problem \eqref{eqn dirchlet prob cross-sec} on the cross-section $\omega_2$.

\begin{Lem}\label{Lemma1}
There exists a positive constant $C=C(p,a,b,|\omega_1|)$ such that
\begin{equation}
\norm{\D u_\ell}_{L^p(\Omega_\ell)}\leq C\ell^{\frac{m}{p}}\norm{f}_{L^{p'}(\omega_2)}^{\frac{1}{p-1}}.
\end{equation}
\end{Lem}
\begin{proof}
    Choosing $u_{\ell}$ as the test function in the weak formulation satisfied by $u_\ell$ itself, using Lemma \ref{lemma_properties of H}, the H\"older and Poincaré inequalities, we get
\begin{equation*}
\begin{split}
    \theta_{1}^{p}\norm{\D \ul}^{p}_{L^{p}(\om_{\ell})}&\leq \int_{\om_{\ell}}H^{p-1}(\D \ul) \left<\D H(\D u_{\ell})\cdot\D \ul\right> \ dx \\
    &= \int_{\om_{\ell}}f(X_{2})\ul\ dx\\
    &\leq \norm{\ul}_{L^{p}(\om_{\ell})}\left(\int_{\om_{\ell}}|f(X_{2})|^{\frac{p}{p-1}} \ dX_{2} \right)^{1-\frac{1}{p}}\\
    &\leq C_{1}\norm{\D \ul}_{L^{p}(\om_{\ell})}\ell^{\frac{m(p-1)}{p}}|\omega_{1}|^{\frac{p-1}{p}}\norm{f}_{L^{\frac{p}{p-1}}(\2)},
\end{split}
\end{equation*}
where $C_{1}$ is the Poincar\'e constant on $\Omega_l$, see Lemma \ref{poincare on strip}. Thus, with $C=\frac{C_{1}|\1|}{a^{p}}$ we have
\begin{equation*}
    \norm{\D \ul}_{L^{p}(\om_{\ell})}\leq C \ell^{\frac{m}{p}}\norm{f}_{L^{p'}(\2)}^{\frac{1}{p-1}}.
\end{equation*}
This finishes the proof of the lemma.
\end{proof}

The proof of the following lemma relies on the same strategy as that used in \cite{chipotpaper2024}. We recall it here for the sake of completeness.
\begin{Lem}\label{Lemma2}
Let $1<p<2$ and $\alpha>\frac{4}{2-p}$. Then there exists a positive constant $C=C(p,c_3)$ such that
\begin{multline*}
\int_{\Omega_\ell}\left|H^{p-1}(\D u_\ell)\D H(\D u_\ell)-H^{p-1}(\D_{2} u_\infty)\D H(\D_{2} u_\infty)\right|\left|u_\ell-u_\infty\right|\rho^{\alpha-1}dx \\ \leq C I_p^{\frac{1}{2}+\frac{1}{p'}}\left(\int_{\Omega_\ell} (|\D u_\ell|+|\D u_\infty|)^p\,dx\right)^{\frac{1}{p}-\frac{1}{2}},
\end{multline*}
where $\rho\in C^1_c(\Omega_\ell)$ is such that $0\leq \rho\leq 1$ and 
\begin{equation*}
I_p=\int_{\Omega_\ell}\left(\left|\D u_\ell\right|+\left|\D u_\infty\right|\right)^{p-2}\left|\D \left(u_\ell- u_\infty\right)\right|^2\rho^\alpha\,dx.
\end{equation*}
\end{Lem}
\begin{proof}
Using \eqref{a2} and H\"older's inequality, we infer that
\begin{align*}
& \int_{\Omega_\ell}\left|H^{p-1}\left(\D u_\ell\right)\D H\left(\D u_\ell\right)-H^{p-1}\left(\D_{2} u_\infty\right)\D H\left(\D_{2} u_\infty\right)\right|\left|u_\ell-u_\infty\right|\rho^{\alpha-1}\,dx\\
&\leq C\int_{\Omega_\ell}\left|\D u_\ell-\D u_\infty\right|\left(\left|\D u_\ell\right|+\left|\D u_\infty\right|\right)^{p-2}|u_\ell-u_\infty|\rho^{\frac{\alpha}{p'}}\rho^{\frac{\alpha}{p}-1}\,dx\\
& \leq C \left(\int_{\Omega_\ell}\left|u_\ell-u_\infty\right|^p\rho^{\alpha-p}\,dx\right)^{\frac{1}{p}}\bigg(\int_{\Omega_\ell}(|\D u_\ell|+|\D u_\infty|)^{(p-2)p'}|\D(u_\ell-u_\infty)|^{p'}\rho^\alpha\,dx\bigg)^{\frac{1}{p'}}\\
& \leq C I_{1}I_{2}.
\end{align*}
To handle $I_{2}$, we estimate using the triangle inequality $|\D\ul-\D\ui|^{p'}\leq |\D\ul-\D\ui|^{2}\left(|\D \ul|+|\D\ui|\right)^{p'-2}$, and observing that $(p-2)p'+p'-2=p-2$, we obtain
\begin{equation*}
    I_2\leq I_p^{\frac{1}{p'}}.
\end{equation*}
%On one hand, we have:
%\begin{align}\label{eq4}
%& \bigg(\int_{\Omega_\ell}(|\D u_\ell|+|\D u_\infty|)^{(p-2)p'}|\D(u_\ell-u_\infty)|^{p'}\rho^\alpha\,dx\bigg)^{\frac{1}{p'}}\notag\\
%& = \bigg(\int_{\Omega_\ell}(|\D u_\ell|+|\D u_\infty|)^{(p-2)p'}|\D(u_\ell-u_\infty)|^2|\D(u_\ell-u_\infty)|^{p'-2}\rho^\alpha\,dx\bigg)^{\frac{1}{p'}}\leq I_p^{\frac{1}{p'}}.
%\end{align}
For $I_{1}$, using Poincaré's inequality on strips, H\"older's inequality and the fact that $\frac{2}{p}(\alpha-p)\geq\alpha$, we deduce that
\begin{align}\label{inequality in lemma 2}
& \left(\int_{\Omega_\ell}|u_\ell-u_\infty|^p\rho^{\alpha-p}\,dx\right)^{\frac{1}{p}}\leq \tilde{C}_{p} \left(\int_{\Omega_\ell}\left|\D(u_\ell-u_\infty)\right|^p\rho^{\alpha-p}\,dx\right)^{\frac{1}{p}}\\
& =C \left(\int_{\Omega_\ell}\left(|\D u_\ell|+|\D u_\infty|\right)^{\frac{(p-2)p}{2}}\left|\D (u_\ell-u_\infty)\right|^p\rho^{\alpha-p}\left(|\D u_\ell|+|\D u_\infty|\right)^{\frac{(2-p)p}{2}}\,dx\right)^{\frac{1}{p}} \notag\\
& \leq C \left(\int_{\Omega_\ell}\left(|\D u_\ell|+|\D u_\infty|\right)^{p-2}\left|\D (u_\ell-u_\infty)\right|^2\rho^{\frac{2}{p}(\alpha-p)}\,dx\right)^{\frac{1}{2}}\left(\int_{\Omega_\ell}\left(|\D u_\ell|+|\D u_\infty|\right)^p\,dx\right)^{\frac{1}{p}-\frac{1}{2}}\notag\\
& \leq C I_p^\frac{1}{2}\left(\int_{\Omega_\ell}\left(|\D u_\ell|+|\D u_\infty|\right)^p\,dx\right)^{\frac{1}{p}-\frac{1}{2}}.\notag
\end{align}
Combining the inequalities for $I_{1}$ and $I_{2}$, we obtain the thesis.
\end{proof}
We next prove a lemma, which is useful in proving Theorem \ref{thm expo convergence}.
\begin{Lem}\label{Lemma 3}
    For any $\phi\in \W(\om_{\ell})$, one has
    \begin{equation}
\label{e2}
\int_{\Omega_\ell}\left<\big(H^{p-1}(\D u_\ell)\D H(\D u_\ell)-H^{p-1}(\D u_\infty)\D H(\D u_\infty)\big)\cdot \D \phi\right>\,dx=0.
\end{equation}
\end{Lem}
\begin{proof}
    On one hand, $u_\infty$ satisfies the weak formulation
\begin{equation}\label{eq1}
\int_{\omega_2}H^{p-1}(0,\D_{X_2}u_\infty)\left<\D_{Z_2}H(0,\D_{X_2} u_\infty)\cdot\D_{X_2}\phi\right>\,dX_2=\int_{\omega_2}f(X_2)\phi\,dX_2,
\end{equation}
for all $\phi\in W^{1,p}_0(\Omega_\ell)$ and for a.e. $X_1\in\omega_1$.
Since $u_\infty$ is constant with respect to the variables of $X_1$, by Fubini's theorem, we derive that 
\begin{align}\label{eq2}
& \int_{\Omega_\ell}H^{p-1}(\D u_\infty)\left <\D_{Z_1} H(\D u_\infty)\cdot\D_{X_1} \phi\right>\,dx\\
& = \sum_{i=1}^m\int_{\omega_2}H^{p-1}(\D u_\infty)\D_{z_i} H(\D u_\infty)\left(\int_{\omega_1}\D_{x_i} \phi\,dX_1\right)\,dX_2=0. \nonumber
\end{align}
Thus, combining \eqref{eq1} and \eqref{eq2}, we get
\begin{align*}
\int_{\Omega_\ell}H^{p-1}(\D u_\infty)\left<\D H(\D u_\infty)\cdot\D \phi\right>\,dx=\int_{\Omega_\ell}f(X_2)\phi\,dx,\quad\forall\phi\in W^{1,p}_0(\Omega_\ell).
\end{align*}
On the other hand, $u_\ell$ satisfies the weak formulation
\begin{equation*}
\int_{\Omega_\ell}H^{p-1}(\D u_\ell)\left<\D H(\D u_\ell)\cdot\D \phi\right>\,dx=\int_{\Omega_\ell}f(X_2)\phi\,dx,\quad\forall\phi\in W^{1,p}_0(\Omega_\ell).
\end{equation*}
The lemma follows by subtracting the previous two equations.
\end{proof}
\subsection{Proof of the Theorem \ref{thm poly convergence}}
\begin{proof}
We divide the proof into two cases; the first is for $p\geq 2$, and the second is for $1<p<2$. The first case depends on Lemma \ref{Lemma1}, whereas the latter follows from Lemma \ref{Lemma1} and Lemma \ref{Lemma2}. \smallskip

Let $\alpha>1$. We choose $\phi(x):=\rho_\ell^{\alpha}(X_1)(u_\ell(x)-u_\infty(x))$ in \eqref{e2}, for $x\in\Omega_\ell$, where $\rho_\ell\in C_c^\infty(\ell\omega_1)$ is a cut-off function such that $\rho_\ell=1$ on $(\ell/2)\omega_1$ and $|\D_{X_1}\rho_\ell|\leq 2/\ell$. It follows that
\begin{align}
\label{eq6}
&\int_{\Omega_\ell}\rho_\ell^\alpha\left<\big(H^{p-1}(\D u_\ell)\D H(\D u_\ell)-H^{p-1}(\D u_\infty)\D H(\D u_\infty)\big)\cdot\D(u_\ell-u_\infty)\right>\,dx\notag\\
& =-\alpha\int_{\Omega_\ell}\rho_\ell^{\alpha-1}\left<\big(H^{p-1}(\D u_\ell)\D_{Z_{1}} H(\D u_\ell)-H^{p-1}(\D u_\infty)\D_{Z_{1}} H(\D u_\infty)\big)\cdot(u_\ell-u_\infty)\D_{X_{1}}\rho_\ell\right>\,dx.
\end{align}
We need to distinguish two cases.\\
\indent \textbf{Case 1: $p\geq 2$.} We set $\alpha:=\frac{p}{p-1}$. Using \eqref{monotonicity}, by \eqref{eq6}, Lemma \ref{inequlity lemma} and Poincaré's inequality on strips (see Lemma \ref{poincare on strip})., we have
\begin{align}
\label{eq3}
&\int_{\Omega_\ell}|\D (u_\ell-u_\infty)|^p\rho_\ell^{\frac{p}{p-1}}\,dx\notag\\
& \leq C \int_{\Omega_\ell}\rho_\ell^{\frac{1}{p-1}}  |H^{p-1}(\D u_\ell)\D_{Z_1}H(\D u_\ell)+H^{p-1}(\D u_\infty)\D_{Z_1}H(\D u_\infty)||\D_{X_1}\rho_\ell||u_\ell-u_\infty|\,dx\notag\\
& \leq \frac{C}{\ell} \int_{\Omega_\ell}\rho_\ell^{\frac{1}{p-1}} \big(|\D u_\ell|^{p-1}+|\D u_\infty|^{p-1}\big)|u_\ell-u_\infty|\,dx\notag\\
& \leq \frac{C}{\ell} \left(\int_{\Omega_\ell}|u_\ell-u_\infty|^p\rho_\ell^{\frac{p}{p-1}}\,dx\right)^{\frac{1}{p}}\left(\int_{\Omega_\ell}\left(|\D u_\ell|^p+|\D u_\infty|^p\right)\,dx\right)^{\frac{p-1}{p}}\notag \\
& \leq \frac{C}{\ell} \left(\int_{\Omega_\ell}|\D_{X_2}(u_\ell-u_\infty)|^p\rho_\ell^{\frac{p}{p-1}}\,dx\right)^{\frac{1}{p}}\left(\int_{\Omega_\ell}\left(|\D u_\ell|^p+|\D u_\infty|^p\right)\,dx\right)^{\frac{p-1}{p}}\notag\\
& \leq \frac{C}{\ell}\left(\int_{\Omega_\ell}|\D(u_\ell-u_\infty)|^p\rho_\ell^{\frac{p}{p-1}}\,dx\right)^{\frac{1}{p}}\bigg(\int_{\Omega_\ell}\big(|\D u_\ell|^p+|\D u_\infty|^p\big)\,dx\bigg)^{\frac{p-1}{p}}.
\end{align}
 Applying Lemma \ref{Lemma1}, we infer that
%\begin{align*}
%\bigg(\int_{\Omega_\ell}|\D u_\ell|^p\,dx\bigg)^{\frac{p-1}{p}}\leq c(p,a,b,c_1,|\omega_1|)\ell^{\frac{m(p-1)}{p}}\bigg(\int_{\omega_2}|\D_{X_2}u_\infty|^p\,dX_2\bigg)^{\frac{p-1}{p}}.
%\end{align*}
%Inserting the previous estimate in \eqref{eq3}, since $u_\infty$ is constant with respect to the variables $X_1$, we get

\begin{align*}
\int_{\Omega_\ell}|\D (u_\ell-u_\infty)|^p\rho_\ell^{\frac{p}{p-1}}\,dx
 \leq \frac{C}{\ell}\left(\int_{\Omega_\ell}|\D(u_\ell-u_\infty)|^p\rho_\ell^{\frac{p}{p-1}}\,dx\right)^{\frac{1}{p}}\ell^{\frac{m(p-1)}{p}}\Bigg(\norm{f}_{L^{p'}(\omega_2)}+\norm{\D_{X_2}u_\infty}^{p-1}_{L^{p}(\omega_2)}\Bigg).
\end{align*}
Finally, by the properties of $\rho_\ell$, we get

\begin{equation}
\int_{\Omega_\frac{\ell}{2}}|\D (u_\ell-u_\infty)|^p\,dx\leq C \ell^{m-\frac{p}{p-1}}\left(\norm{f}_{L^{p'}(\omega_2)}^{p'}+\norm{\D_{X_2}u_\infty}_{L^p(\omega_2)}^p\right),
\end{equation}

which concludes the proof in the first case.\smallskip

\indent\textbf{Case 2: $1<p<2$.} Let $\alpha>\frac{4}{2-p}$. By assumption \eqref{a1}, \eqref{eq6}, Lemma \ref{Lemma2} and Young's inequality with $\epsilon>0$, it follows that
\begin{align*}
& I_p=\int_{\Omega_\ell}(|\D u_\ell|+|\D u_\infty|)^{p-2}|\D (u_\ell-u_\infty)|^2\rho_\ell^\alpha\,dx\\
%& \leq c(c_2)\int_{\Omega_\ell}\rho_\ell^{\alpha}\prodscal{H^{p-1}(\D u_\ell)\D H(\D u_\ell)-H^{p-1}(\D u_\infty)\D H(\D u_\infty)}{\D(u_\ell-u_\infty)}\,dx\\
%& =-c(p,c_2)\int_{\Omega_\ell}\rho_\ell^{\alpha-1}\prodscal{H^{p-1}(\D u_\ell)\D H(\D u_\ell)-H^{p-1}(\D u_\infty)\D H(\D u_\infty)}{(u_\ell-u_\infty)\D\rho_\ell}\,dx\\
& \leq \frac{C}{\ell} \int_{\Omega_\ell}|H^{p-1}(\D u_\ell)\D H(\D u_\ell)-H^{p-1}(\D u_\infty)\D H(\D u_\infty)||u_\ell-u_\infty|\rho^{\alpha-1}\,dx\\
& \leq \frac{C}{\ell} I_p^{\frac{1}{2}+\frac{1}{p'}}\left(\int_{\Omega_\ell} \left(|\D u_\ell|+|\D u_\infty|\right)^p\,dx\right)^{\frac{1}{p}-\frac{1}{2}}\\
& \leq C\left(\varepsilon I_p+\frac{1}{\varepsilon \ell^{\frac{2p}{2-p}}}\int_{\Omega_\ell} \left(|\D u_\ell|^p+|\D u_\infty|^p\right)\,dx\right),
\end{align*}
for any $\varepsilon>0$. Choosing $\varepsilon$ sufficiently small and applying Lemma \ref{Lemma1}, since $u_\infty$ and $f$ does not depend on $X_1$, we deduce that
\begin{align*}
I_p\leq \frac{C}{\ell^{\frac{2p}{2-p}}}\int_{\Omega_\ell} \big(|\D u_\ell|^p+|\D u_\infty|^p\big)\,dx\leq \frac{C}{\ell^{\frac{2p}{2-p}}}\ell^m\bigg(\norm{f}^{p'}_{L^{p'}(\omega_2)}+\norm{\D_{X_2}u_\infty}^p\bigg).
\end{align*}
Thus, reasoning as in \eqref{inequality in lemma 2} and taking into account the choice of $\alpha$ and the previous chain of inequalities, we finally get that
\begin{align*}
& \int_{\Omega_{\ell/2}}|\D(u_\ell-u_\infty)|^p\,dx\\
& \leq\int_{\Omega_\ell}|\D(u_\ell-u_\infty)|^p\rho^{\alpha-p}\,dx\\
& 
=\int_{\Omega_\ell}(|\D u_\ell|+|\D u_\infty|)^{\frac{(p-2)p}{2}}|\D(u_\ell-u_\infty)|^p\rho^{\alpha-p}(|\D u_\ell|+|\D u_\infty|)^{\frac{(2-p)p}{2}}\,dx\\
& \leq \Bigg(\int_{\Omega_\ell}(|\D u_\ell|+|\D u_\infty|)^{p-2}|\D(u_\ell-u_\infty)|^2\rho^{\frac{2}{p}(\alpha-p)}\,dx\Bigg)^{\frac{p}{2}}\Bigg(\int_{\Omega_\ell}(|\D u_\ell|+|\D u_\infty|)^p\,dx\Bigg)^{1-\frac{p}{2}}\\
& \leq CI_p^\frac{p}{2}\left(\int_{\Omega_\ell}(|\D u_\ell|+|\D u_\infty|)^p\,dx\right)^{1-\frac{p}{2}}\\
& \leq \frac{C}{\ell^{\frac{p^2}{2-p}}}\int_{\Omega_\ell} \big(|\D u_\ell|^p+|\D u_\infty|^p\big)\,dx\leq \ell^{m-\frac{p^2}{2-p}}\bigg(\norm{f}_{L^p(\omega_2)}^{p'}+\norm{\D_{X_2}u_\infty}_{L^p(\omega_2)}^p\bigg),
\end{align*}

which leads to the thesis with $C=C\left(p,\theta_{1},\theta_{2},c_j,\left|\omega_1\right|,|\2|,\norm{\D_{X_{2}}\ui}_{L^{p}(\2)}^{p},\norm{f}_{L^p(\omega_2)}^{p'}\right)$.
\end{proof}

\subsection{Exponential rate of convergence}

For a fixed $p\geq 2$ and any $1<q<+\infty$ let us consider 
\begin{equation}\label{Eqn restricted class of H}
\Tilde{H}_{q}\left(z\right)=\left(F^{q}\left(Z_1\right)+G^{q}\left(Z_2\right)\right)^{\frac{1}{q}},
\end{equation} with $z=\left(Z_{1},Z_{2}\right)$, where, $F:\R^{r}\to \R$ and $G:\R^{n-r}\to \R$ satisfy
\begin{enumerate}
    \myitem{$\textbf{A3}$} \label{a3} \textit{Monotonicity
    \begin{align*}
&\left<F^{p-1}\left(Z_1\right)\D F\left(Z_1\right)-F^{p-1}\left(W_1\right)\D F\left(W_1\right) \cdot Z_1-W_1\right>\geq c_1\left|Z_1-W_1\right|^p,\\
&\left<G^{p-1}\left(Z_2\right)\D G\left(Z_2\right)-G^{p-1}\left(W_2\right)\D G\left(W_2\right)\cdot Z_2-W_2\right>\geq c_1\left|Z_2-W_2\right|^p,
\end{align*}
and 
\begin{align*}
&\theta_{1}\left|Z_1\right|\leq F\left(Z_1\right)\leq \theta_{2}\left|Z_1\right|,\\
&\theta_{1}\left|Z_2\right|\leq G\left(Z_2\right)\leq \theta_{2}\left|Z_2\right|,
\end{align*}
for all $z\in\R^n$}.
\end{enumerate}
Then, in the case that $q=p$, we have the following result.
\begin{Thm}\label{thm expo convergence}
Let $p\geq 2$ and $\ul,\ui$ solve \eqref{eqn dirichlet problem on omel}, \eqref{eqn dirchlet prob cross-sec} respectively with $H=\Tilde{H}_{p}$ as in \eqref{Eqn restricted class of H}, then 
\begin{equation*}
\label{H}
\norm{\D(u_\ell-u_\infty)}_{L^p\left(\Omega_{\ell/2}\right)}\leq C e^{-\beta \ell},
\end{equation*}
for some positive constants $C$ and $\beta$ both independent of $\ell$.
\end{Thm}
The following proof is based on the well-known \textit{hole-filling technique}.
\begin{proof}
We denote by $\ell_1$ a real number such that $0<\ell_1<\ell-1$. For $\alpha>0$, we set $\phi(x):=\rho_{\ell_1}^{\alpha}(X_1)(u_\ell(x)-u_\infty(x))$ in \eqref{e2}, for $x\in\Omega_\ell$, where $\rho_{\ell_1}\in C_c^\infty((\ell_1+1)\omega_1)$ is a cut-off function such that $\rho_{\ell_1}=1$ on $\ell_1\omega_1$ and $|\D_{X_1}\rho_\ell|\leq c$, for some positive constant $c=c(\omega_1)$. It holds that
\begin{equation}
H^{p-1}(z)\D_{Z_1} H(z)=F^{p-1}(Z_1)\D_{Z_1} F(Z_1),\quad\forall z\in\R^n.
\end{equation}
By the properties of $\rho_{\ell_1}$, \eqref{a1}, \eqref{eq6} and since $F(\D_{X_{1}}\ui)=0$ as $\D_{X_{1}}\ui=0$, applying  Young's inequality and the Poincaré's inequality on strips, we infer that
\begin{align*}
& \int_{\Omega_{\ell_1}}\left|\D (u_\ell-u_\infty)\right|^p\,dx\leq\int_{\Omega_{\ell_1+1}}|\D (u_\ell-u_\infty)|^p\rho_{\ell_1}^{\alpha}\,dx\\
& \leq C\int_{\Omega_{\ell_1+1}}\rho_{\ell_1}^{\alpha}\left<{H^{p-1}(\D u_\ell)\D H(\D u_\ell)-H^{p-1}(\D u_\infty)\D H(\D u_\infty)}\cdot{\D(u_\ell-u_\infty)}\right>\,dx\\
& =-\alpha C\int_{\Omega_{\ell_1+1}\setminus\Omega_{\ell_1}}\rho_{\ell_1}^{\alpha-1}\left<{H^{p-1}(\D u_\ell)\D_{Z_1} H(\D u_\ell)-H^{p-1}(\D u_\infty)\D_{Z_1} H(\D u_\infty)\cdot(u_\ell-u_\infty)\D_{X_1}\rho_\ell}\right>\,dx\\
%\setminus\Omega_{\ell_1}}
& \leq C\int_{\Omega_{\ell_1+1}\setminus \Omega_{\ell_{1}}}\rho_{\ell_1}^{\alpha-1}\left|\left<F^{p-1}\left(\D_{X_1} u_\ell\right)\D_{Z_1} F\left(\D_{X_1} u_\ell\right) \cdot\left(u_\ell-u_\infty\right)\D_{X_1}\rho_\ell\right>\right|\\
& +\left|\left<F^{p-1}\left(\D_{X_1} u_\infty\right)\D_{Z_1} F\left(\D_{X_1} u_\infty\right) \cdot\left(u_\ell-u_\infty\right)\D_{X_1}\rho_\ell\right>\right|\, dx\\
& \leq C\int_{\Omega_{\ell_1+1}\setminus\Omega_{\ell_1}}|\D_{X_1} u_\ell|^{p-1}|u_\ell-u_\infty|\,dx\\
& =C\int_{\Omega_{\ell_1+1}\setminus\Omega_{\ell_1}}|\D_{X_1} (u_\ell-u_\infty)|^{p-1}|u_\ell-u_\infty|\,dx\\
& \leq C \left(\int_{\Omega_{\ell_1+1}\setminus\Omega_{\ell_1}}|\D_{X_1} (u_\ell-u_\infty)|^p\,dx+\int_{\Omega_{\ell_1+1}\setminus\Omega_{\ell_1}}|u_\ell-u_\infty|^p\,dx\right)\\
& \leq C \int_{\Omega_{\ell_1+1}\setminus\Omega_{\ell_1}}|\D (u_\ell-u_\infty)|^p\,dx.
\end{align*}
%\begin{align*}
%& \int_{\Omega_{\ell_1}}|\D (u_\ell-u_\infty)|^p\,dx\\
%& \leq c(\alpha,b,p,c_1,\omega_1)\bigg(\int_{\Omega_{\ell_1+1}\setminus\Omega_{\ell_1}}|\D_{X_1} (u_\ell-u_\infty)|^p\,dx+\int_{\Omega_{\ell_1+1}\setminus\Omega_{\ell_1}}|u_\ell-u_\infty|^p\,dx\bigg)\\
%& \leq c(\alpha,b,p,c_1,\omega_1,\omega_2)\bigg(\int_{\Omega_{\ell_1+1}\setminus\Omega_{\ell_1}}|\D_{X_1} (u_\ell-u_\infty)|^p\,dx+\int_{\Omega_{\ell_1+1}\setminus\Omega_{\ell_1}}|\D_{X_2}(u_\ell-u_\infty)|^p\,dx\bigg)\\
%& \leq c(\alpha,b,p,c_1,\omega_1,\omega_2)\int_{\Omega_{\ell_1+1}\setminus\Omega_{\ell_1}}|\D (u_\ell-u_\infty)|^p\,dx.
%\end{align*}
By ``filling the hole" on the right side and dividing by $C+1$  both sides, we obtain
\begin{equation*}
\int_{\Omega_{\ell_1}}|\D (u_\ell-u_\infty)|^p\,dx\leq \frac{C}{C+1}\int_{\Omega_{\ell_1+1}}|\D (u_\ell-u_\infty)|^p\,dx.
\end{equation*}
Setting $\gamma:=\frac{C}{C+1}$, starting from $\ell_1=\frac{\ell}{2}$ and iterating $[\frac{\ell}{2}]$ times, where $[\cdot]$ denotes the integer part of a real number, we gain,
\begin{equation*}
\int_{\Omega_{\ell/ 2}}|\D (u_\ell-u_\infty)|^p\,dx\leq \gamma^{[\frac{\ell}{2}]}\int_{\Omega_{\frac{\ell}{2}+[\frac{\ell}{2}]}}|\D (u_\ell-u_\infty)|^p\,dx\leq \gamma^{\frac{\ell}{2}-1}\int_{\Omega_{\ell}}|\D (u_\ell-u_\infty)|^p\,dx.
\end{equation*}
By Lemma \ref{Lemma1}, we get
\begin{align*}
\norm{\D (u_\ell-u_\infty)}_{L^p(\Omega_{\ell/2})}
& \leq C\gamma^{\frac{1}{p}\big(\frac{\ell}{2}-1\big)}\ell^{\frac{m}{p}}\left(\norm{f}_{L^{p'}(\omega_2)}^{\frac{1}{p-1}}+\norm{\D_{X_2}u_\infty}_{L^p(\omega_2)}\right)\\
& \leq C e^{\ell\frac{\log{\gamma}}{2p}}\ell^{\frac{m}{p}} \leq Ce^{-\beta\ell},
\end{align*}
for any $\beta\in\Big(0,-\frac{\log\gamma}{2p}\Big)$, which leads to the result.
\end{proof}

\section{Energy asymptotics}\label{Energy asymptotics}
In this section, we deal with the proof of Theorem \ref{thrm energy}, whose proof is quite similar to that given in \cite{PJana2024} for the anisotropic $p$-Laplacian. The solutions $\ul$ and $\ui$ are minimisers of functionals $J_{\ell}$ and $J_{\infty}$ respectively, that is
\begin{equation*}
    J_{\ell}\left(u_{\ell}\right)=\displaystyle \min_{v\in W^{1,p}_{0}\left(\om_{\ell}\right)} J_{\ell}\left(v\right)  \quad\text{and}\quad J_{\infty}\left(u_{\infty}\right)=\min_{w\in W^{1,p}_{0}\left(\2\right)}J_{\infty}\left(w\right).
\end{equation*}
It follows the proof of Theorem \ref{thrm energy}.
\begin{proof}
    We construct an admissible test function for the left-side inequality using $\ul$ and the Jensen's inequality. The second one is choosing the right test function and some manipulation.
    \smallskip
    
    By taking 
 \begin{equation*}
     w_{\ell}\left(X_{2}\right)=\frac{1}{|\ell \omega_{1}|}\int_{\ell\omega_{1}}u_{\ell}\left(X_{1},X_{2}\right)\ dX_{1},
 \end{equation*}
one has by H\"older inequality that $w_{\ell}\in W_{0}^{1,p}(\omega_{2})$. Thus, by Jensen's inequality, homogeneity and Lemma \ref{lemma_properties of H}, we estimate
\begin{equation*}
    \begin{split}
        J_{\infty}\left(u_{\infty}\right)&\leq J_{\infty}\left(w_{\ell}\right)=\int_{\omega_{2}}\frac{1}{p}H^{p}\left(\D_{2} w_{\ell}\right)-f\left(X_{2}\right)w_{\ell} \ dX_{2} \\
        &= \int_{\omega_{2}}\frac{1}{p}H^{p}\left(\frac{1}{\left|\ell
        \1\right|}\int_{\ell\omega_{1}}\D_{2} u_{\ell}\left(X_{1},X_{2}\right)\ dX_{1}\right)\ dX_{2}-\frac{1}{\left|\ell\1\right|}\int_{\om_{\ell}}f\left(X_{2}\right)u_{\ell}(x) \ dx \\
        &\leq \frac{1}{\left|\ell \omega_{1}\right|}\left(\int_{\om_{\ell}}\frac{1}{p}H^{p}(\D u_{\ell})\ dx -\int_{\om_{\ell}}f(X_{2})u_{\ell}(x)\,dx\right)\\
        &\leq \frac{J_{\ell}\left(u_{\ell}\right)}{\left|\ell\omega_{1}\right|}.
    \end{split}
\end{equation*}
By this, we see that the scaling $\frac{1}{|\ell\omega_{1}|}$ normalizes the missing coordinates from $\omega_{2}$.\smallskip

For the other inequality, choose a function $\phi_{\ell}\in C^{1}(\R^{m})$ such that $0\leq \phi\leq 1$, $\phi=1$ on $(\ell-1)\omega_{1}$, $\phi=0$ in $(\ell \omega_{1})^{c}$ and $|\D \phi|\leq M$, for some positive constant $M$. Since $\phi(X_{1})u_{\infty}(X_{2})\in W_{0}^{1,p}(\om_{\ell})$,
    \begin{equation*}
        \begin{split}
            J_{\ell}(u_{\ell})
            &\leq J_{\ell}\left(\phi u_{\infty}\right)=J_{\ell-1}(u_{\infty})+\int_{\om_{\ell}\setminus\om_{\ell-1}} \frac{1}{p}H^{p}(\D (\phi(X_{1})u_{\infty}(X_{2})))-f\left(X_{2}\right)\phi\left(X_{1}\right) u_{\infty}\left(X_{2}\right) \ dx    \\
            & = J_{\ell}\left(u_{\infty}\right)+\int_{\om_{\ell}\setminus\om_{\ell-1}} \frac{1}{p}\big[H^{p}\left(\D \left(\phi u_{\infty}\right)\right)-H^{p}\left(\left(0,\D_{2}u_{\infty}\right)\right)\big]-fu_{\infty}\left(\phi-1\right) \ dx\\
           &\leq \left|\ell\omega_{1}\right| J_{\infty}\left(u_{\infty}\right)+\int_{\om_{\ell}\setminus\om_{\ell-1}} \frac{1}{p}\big[H^{p}\left(\D \left(\phi u_{\infty}\right)\right)+H^{p}\left(0,\D_{2}u_{\infty}\right)\big]+\left|fu_{\infty}\right| \ dx\\
           & \leq \left|\ell\omega_{1}\right| J_{\infty}\left(u_{\infty}\right)\\
           &\quad \quad+\frac{\theta_2^{p} }{p}\int_{\om_{\ell}\setminus\om_{\ell-1}} \big[|u_{\infty}\D_{1}\phi|^p+|\phi\D_{2}u_{\infty}|^{p}\big]\,dx+ \{|\ell\omega_{1}|-|\left(\ell-1\right)\omega_{1}|\}\norm{f}_{L^{p'}(\omega_{2})}\norm{u_{\infty}}_{L^{p}(\omega_{2})}\\
           %& \leq |\ell\omega_{1}| J_{\infty}(u_{\infty})+\int_{\om_{\ell}\setminus\om_{\ell-1}} \frac{b^{p}}{p} \left|M u_{\infty}+\D_{2}u_{\infty}\right|^{p}\ dx+ \{|\ell\omega_{1}|-|(\ell-1)\omega_{1}|\}\norm{f}_{L^{p'}(\omega_{2})}\norm{u_{\infty}}_{L^{p}(\omega_{2})}\\
           & \leq |\ell\omega_{1}| J_{\infty}(u_{\infty})+\{|\ell\omega_{1}|-|(\ell-1)\omega_{1}|\} \frac{\theta_2^{p}}{p}\int_{\omega_{2}} \big[|M \ui|^p+|\D_{2}u_{\infty}|^{p}\big]\,dx + \norm{f}_{L^{p'}(\omega_{2})}\norm{u_{\infty}}_{L^{p}(\omega_{2})}\\
           & \leq |\ell\omega_{1}| J_{\infty}(u_{\infty})+\{|\ell\omega_{1}|-|(\ell-1)\omega_{1}|\}C.
       \end{split}
    \end{equation*}
    Thus we have 
    \begin{equation*}
        \frac{J_{\ell}(u_{\ell})}{|\ell\omega_{1}|}\leq J_{\infty}(u_{\infty})+\frac{C}{\ell},
    \end{equation*}
    which completes the proof.
\end{proof}

\section{An eigenvalue problem}\label{An eigenvalue problem}
In this section, we will show that the sequence of eigenvalues $\{\lambda_{\ell}^{1}\}_{\ell}$ converges to the eigenvalue $\mu_{\infty}$ of the cross-sectional problem. 
\subsection{Proof of Theorem \ref{eigenvalue from above}}
\begin{proof}
     First, we establish the upper bound. Let $\phi:\ell\omega_{1}\to \R$ be a Lipschitz continuous function such that $0\leq \phi\leq 1$, $\phi=0$ on $\partial (\ell\omega_{1})$, $\phi=1$ on $\frac{\ell}{2}\omega_1$ and $|\D_{1} \phi|\leq \frac{2}{\ell}$. Then $\phi u_{\infty}\in W_{0}^{1,p}(\om_{\ell})$, thus, using the sub-linearity of $H$, taking $c=1+\frac{1}{\ell}$ in (2) of Lemma \ref{inequlity lemma}, one can derive
    \begin{equation}
    \label{eqq1}
        \begin{split}
            \lambda_{\ell}^{1} \int_{\om_{\ell}}\left|\phi u_{\infty}\right|^{p} \ dx &\leq  \int_{\om_{\ell}}H^{p}\left(\D\left(\phi u_{\infty}\right)\right) \ dx\\
            &\leq \int_{\om_{\ell}}\left(H\left(\phi \D_{2}u_{\infty}\right)+H\left(u_{\infty}\D_{1}\phi_\ell\right)\right)^{p} \ dx\\
            &\leq \int_{\om_{\ell}}c H^{p}\left(\phi\D_{2}u_{\infty}\right)+\left(1-c^{\frac{-1}{{p}-1}}\right)^{1-{p}}H^{p}\left(u_{\infty}\D_{1}\phi\right) \ dx\\
            &\leq c\int_{\om_{\ell}}\phi^{p} H^{p}\left(\D_{2}u_{\infty}\right)\,dx+\left(1-c^{\frac{-1}{{p}-1}}\right)^{1-{p}}\int_{\Omega_\ell}u_{\infty}^{{p}}H^{p}\left(\D_{1}\phi\right) \ dx.\\
            & = (I)+(II).
        \end{split}
    \end{equation}
    Concerning the first term, we have
    \begin{equation}
    \label{eqq2}
        (I)= \bigg(1+\frac{1}{\ell}\bigg)\int_{\ell\omega_1} \phi\,dX_1\int_{\omega_2} H^p(\D_2 u_\infty)\,dX_2.
    \end{equation}
    For the second addend, by applying Bernoulli's inequality, we get, assuming $l>1$ in such a way that $c=1+1/\ell<2$,
\begin{equation}
    \frac{c}{\Big(c^{\frac{1}{p-1}}-1\Big)^{p-1}}\leq \frac{2}{\bigg(\Big(1+\frac{1}{\ell}\Big)^{\frac{1}{p-1}}-1\bigg)^{p-1}}\leq 2(p-1)^{p-1}\ell^{p-1}.
    \end{equation}
Thus, we estimate
\begin{align}\label{eqq3}
    (II)\notag
    \leq 2^{p+1}(p-1)^{p-1}\ell^{p-1}\frac{\theta_2^p}{\ell^p}\ell^m|\omega_1|\int_{\omega_2} u_{\infty}^p\,dX_2
    & =2^{p+1}(p-1)^{p-1}2^m\frac{\theta_2^p}{\ell}\int_{\frac{\ell}{2}\omega_1}\phi^p\,dX_1\int_{\omega_2} u_{\infty}^p\,dX_2\notag\\
    & \leq  2^{p+1}(p-1)^{p-1} 2^m\frac{\theta_2^p}{\ell}\int_{\ell\omega_1}\phi^p\,dX_1\int_{\omega_2} u_{\infty}^p\,dX_2.
\end{align}
    Inserting \eqref{eqq2} and \eqref{eqq3} in \eqref{eqq1}, dividing by $\int_{\Omega_\ell}|\phi u_\infty|^p\,dx$, we obtain the desired estimate of the upper bound.\\
    \indent 
Next, we establish the lower bound. Setting with an abuse of notation $\ui\left(X_{1},X_{2}\right)=\ui\left(X_{2}\right)$,  we have 
\begin{equation}\label{u_inf on Omega_ell}
\begin{split}
    \int_{\om_{\ell}}H^{{p-1}}(\D \ui)\left<\D H(\D\ui)\cdot \D \psi\right> \ dx &=\int_{\omega_{2}} H^{{p-1}}(\D_{2}\ui) \left<\D_{Z_{2}}H(\D_{2}\ui) \cdot \int_{\ell\omega_{1}} \D \psi\ dX_{1} \right> \ dX_{2}\\
    &= \mu_{\infty}\int_{\om_{\ell}} \ui^{p-1} \psi \ dx,
    \end{split}
\end{equation}
for any $\psi\in C_{c}^{\infty}(\om_{\ell})$. Let $\phi\in C_{c}^{\infty}\left(\om_{\ell}\right)$. Since $\ui >0$ and $\ui\in C^{1}$, it holds that $\phi^{p}/\ui^{p-1}\in W_{0}^{1,p}(\om_{\ell})$ as $\ui>\alpha(\phi)$ on $\mathrm{supp}(\phi)$, for some $\alpha>0$.  By Lemma \ref{picon}, we infer that 
    \begin{equation*}
        H^{p}\left(\D \phi\right)- \left<H^{p-1}\left(\D \ui\right)\D H\left(\D \ui\right)\cdot\D \left(\frac{\phi^{p}}{\ui^{p-1}}\right)\right>\geq 0.
    \end{equation*}
    By integrating and using \eqref{u_inf on Omega_ell},
    \begin{equation*}
            \int_{\om_{\ell}} H^{p}\left(\D \phi\right)\ dx\geq \mu_{\infty}\int_{\om_{\ell}} \left|\phi\right|^{p}\ dx. 
        \end{equation*}
        By the density, minimising over $W_{0}^{1,p}(\om_{\ell})$ we get the result.
\end{proof}
\subsection{Exact rate of convergence}
In Theorem \ref{eigenvalue from above}, if $H$ is such that $H^{p}\left(Z_{1},Z_{2}\right)\geq H^{p}\left(Z_{1},0\right)+H^{p}\left(0,Z_{2}\right)$ one can derive a better result.
\begin{Thm}\label{eigenvalue lower bound withspeed}
    For $H=\tilde{H}_{q}$ as in \eqref{Eqn restricted class of H} with $r=m$, we infer, for $q\leq p$ that
    \begin{equation*}
         \frac{C_{1}}{\ell^{p}}+\mu_{\infty}\leq \lambda_{\ell}^{1}.
    \end{equation*}
    and for $q\geq p$, that 
    \begin{equation*}
        \lambda_{\ell}^{1}\leq \muinf+ \frac{C_{2}}{\ell^{p}},
    \end{equation*}
    where $C_{1}=C_{1}\left(\theta_1,p\right)$ and $C_{2}=C_{2}\left(\theta_2,m,p\right)$ are positive constants.
\end{Thm}
\begin{proof}
    For $q\leq p$, using the structure of $H$, (3) of Lemma \ref{inequlity lemma} and Poincar\'e inequality on strips and \eqref{raleigh quotient for mu}, we get
    \begin{equation*}
        \begin{split}
            \lamel\int_{\om_{\ell}}u_{\ell}^{p} \ dx &=\int_{\om_{\ell}}H^{p}\left(\D u_{\ell}\right)\ dx\\
            &\geq \int_{\om_{\ell}} F^{p}\left(\D_{X_{1}}\ul\right)+G^{p}\left(\D_{X_{2}}u_{\ell}\right) \ dx\\
            &\geq \frac{C}{\ell^{p}}\int_{\omega_{2}}\int_{\ell\omega_{1}}u_{\ell}^{p}\ dX_{1}dX_{2}+  \int_{\ell\omega_{1}} \frac{\int_{\omega_{2}}G^{p}\left(\D_{2}u_{\ell}\right)dX_{2}}{\int_{\omega_{2}} u_{\ell}^{p}\ dX_{2}} \left(\int_{\omega_{2}} u_{\ell}^{p}\ dX_{2}\right) \ dX_{1}\\
            &\geq \frac{C}{\ell^{p}} \int_{\om_{\ell}}u_{\ell}^{p} \ dx + \mu_{\infty}\int_{\om_{\ell}} u_{\ell}^{p}\ dx,
        \end{split}
    \end{equation*}
    from which the inequality follows. For $q \geq p$, we choose a test function of the form $\phi_{\ell}\ui$ where $\phi_{\ell}:\ell\1\to \R$ is a Lipschitz continuous function with $0\leq \phi_\ell\leq 1$, $\phi_\ell=1$ in $(\ell-1)\1$ and $|\D \phi|\leq \frac{1}{\ell}$. Then we have, by \eqref{raliegh quotient for lamda} and (4) of Lemma \ref{inequlity lemma}, that 
    \begin{equation*}
        \begin{split}
            \lamel \int_{\om_{\ell}} \phi_{\ell}^{p}\ui^{p}\ dx &\leq \int_{\om_{\ell}}H^{p} \left(\ui \D_{X_{1}}\phi_{\ell}+\phi_{\ell} \D_{X{2}}\ui \right) \ dx\\
            & =\int_{\om_{\ell}}\left(F^{q}\left(\ui \D_{X_{1}}\phi_{\ell}\right)+G^{q}\left(\phi_{\ell}\D_{X_{2}}\ui\right)\right)^{\frac{p}{q}}\ dx\\
            &\leq \int_{\om_{\ell}}\ui^{p}F^{p}\left(\D_{X_{1}}\phi\right)+\phi_{\ell}^{p}G^{p}\left(\D_{X_{2}} \ui\right) \ dx\\
            &\leq \int_{\2}\ui^{p}\ dX_{2}\int_{\ell\1}\frac{\theta_2^p}{\ell^{p}}\ dX_{1}+\int_{\ell\1}\phi_\ell^{p}\ dX_{1} \int_{\2}G^{p}\left(\D_{X_{2}}\ui\right) \ dX_{2}.
        \end{split}
    \end{equation*}
    Therefore, we have
    \begin{equation*}
        \begin{split}
                \lamel &\leq \frac{\theta_2^p}{\ell^{p}}\frac{\int_{\ell\1} \ dX_{1}}{\int_{\left(\ell/2\right)\1}}\ dX_{1}+\muinf\\
            &\leq \muinf+ \frac{C_{2}}{\ell^{p}},
        \end{split}
    \end{equation*}
    where $C_{2}=2^m\theta_2^p$.
\end{proof}
\section{Appendix}\label{Appendix}
In this section we briefly deal with the existence and uniqueness results for \eqref{eqn dirichlet problem on omel}. We prove the existence of a solution via the direct method of calculus of variation. The case when $p=2$ can be found in \cite{MezeiVas2019}. The uniqueness follows the assumption \eqref{a1}.\smallskip

The energy of \eqref{eqn dirichlet problem on omel} is 
\begin{equation*}
    J(u)=\int_{\om}\left\{\frac{H^{p}\left(\D u\right)}{p}-fu \right\}\ dx, \quad \text{for all}\ u\in W_{0}^{1,p}(\om).
\end{equation*}
\begin{Lem}
    Let $f\in L^{p'}(\om)$. Then  the functional $J(u)$ if bounded from below. Further, it is coercive.
\end{Lem}
\begin{proof}
    For any $\epsilon>0$, by Young's inequality, Lemma \ref{lemma_properties of H} and Poincaré inequality, one can get
    \begin{equation*}
        \begin{split}
            J(u)&\geq \frac{\theta_1^{p}}{p}\int_{\om} |\D u|^{p}\ dx -\frac{\epsilon}{p}\int_{\om}|u|^{p} \ dx- \frac{1}{p'\epsilon^{p'-1}}\int_{\om}|f|^{p'}\ dx\\
            &\geq \frac{\theta_1^{p}C-\epsilon}{p} \int_{\om}|u|^{p}\ dx - \frac{1}{p'\epsilon^{p'-1}}\int_{\om}|f|^{p'}\ dx.
        \end{split}
    \end{equation*}
    Both the result follows by choosing $\epsilon<\theta_1^{p}C$.
\end{proof}
\begin{Lem}[Sequential lower semicontinuity]
Let $\ul\to u$ in $W_{0}^{1,p}(\om)$, then we have 
\begin{equation*}
    J(u)\leq \displaystyle \liminf_{n\to \infty} J(u_{n}).
\end{equation*}
\end{Lem}
\begin{proof}
    By H\"older's inequality, it holds that
    \begin{equation*}
        \int_{\om} |f||u_{n}-u|\ dx\leq \left(\int_{\om}|f|^{p'}\ dx\right)^{\frac{1}{p'}}\left(\int_{\om}|u_{n}-u|^{p}\ dx\right)^{\frac{1}{p}}\to 0 \quad \text{as}\ n\to\infty.
    \end{equation*}
    Since $H\geq0$ and $H\in C(\R^{n})$, by Fatou's lemma we have 
    \begin{equation*}
        \int_{\om}H^{p}(\D u)\ dx \leq \displaystyle\liminf_{n\to \infty}\int_{\om}H^{p}(\D u_{n})\ dx.
    \end{equation*}
    Thus,
    \begin{equation*}
        \begin{split}
            J(u)&\leq \displaystyle \liminf_{n\to \infty}\int_{\om}\frac{1}{p}H^{p}(\D u_{n}) \ dx - \liminf_{n\to \infty}\int_{\om}fu_{n}\ dx \\
            &\leq \liminf_{n\to\infty}J(u_{n}).
        \end{split}
    \end{equation*}
\end{proof}
The weak lower semicontinuity of $J$ follows from the convexity and lower semicontinuity.

\section*{Acknowledgements}
Luca Esposito and Lorenzo Lamberti are members of the Gruppo Nazionale per l’Analisi Matematica, la Probabilità e le loro Applicazioni (GNAMPA) of the Istituto Nazionale di Alta Matematica (INdAM) and wish to acknowledge financial support from INdAM GNAMPA
Project 2024 “Regolarità per problemi a frontiera libera e disuguaglianze funzionali in
contesto finsleriano”.

 Prosenjit Roy would like to thank the University of Salerno for the warm hospitality where
this work originated.

 N.N.Dattatreya is supported by PMRF grant (2302262).

\renewcommand\refname{Bibliography}
\bibliographystyle{abbrv}
\bibliography{ref}

\end{document}